\documentclass[12pt]{article}
\usepackage{amsmath}
\usepackage{amsfonts}
\usepackage{amssymb}
\usepackage[mathscr]{eucal}

\newif{\ifcomentarios}
\comentariosfalse

\newtheorem{theorem}{Theorem}

\newtheorem{corollary}[theorem]{Corollary}

\newtheorem{lemma}[theorem]{Lemma}

\newtheorem{remark}[theorem]{Remark}

\newtheorem{proposition}[theorem]{Proposition}

\newcommand{\be}{\begin{eqnarray}}
\newcommand{\en}{\end{eqnarray}}
\newcommand{\bee}{\begin{eqnarray*}}
\newcommand{\ene}{\end{eqnarray*}}

\begin{document}

\author{S. L. Carvalho\thanks{%
Supported by FAPESP under grant \#06/60711-4. Email: \texttt{%
silas@if.usp.br}}, \ D. H. U. Marchetti\thanks{%
No thanks. Email: \texttt{marchett@if.usp.br}} \ \& W. F. Wreszinski\thanks{%
Email:\texttt{\ wreszins@fma.if.usp.br}} \\
%EndAName
Instituto de F\'{\i}sica\\
Universidade de S\~{a}o Paulo \\
Caixa Postal 66318\\
05314-970 S\~{a}o Paulo, SP, Brasil}
\title{Pointwise Decay of Fourier-Stieltjes transform of the Spectral
Measure for Jacobi Matrices with Faster-than-Exponential Sparse Perturbations%
}
\date{}
\maketitle

\begin{abstract}
We consider off-diagonal Jacobi matrices $J$ with
(faster--than--exponential) sparse perturbations. We prove (Theorem
\ref{onehalf}) that the Fourier transform $\widehat{\left\vert
f\right\vert 
^{2}d\rho }(t)$ of the spectral measure $\rho $ of $J$, whose sparse
perturbations are at least separated by a distance $\exp \left( cj(\ln
j)^{2}\right) /\delta ^{j}$, for some $c>1/2,$ $0<\delta <1$ and for a dense
subset of $C_{0}^{\infty }(-2,2)$--functions $f$, decays as $t^{-1/2}\Omega
(t)$, uniformly in the spectrum $[-2,2]$, $\Omega (t)$ increasing less
rapidly than any positive power of $t$, improving earlier results obtained
by Simon (Commun. Math. Phys. \textbf{179}, 713-722 (1996)) and by
Krutikov--Remling (Commun. Math. Phys. \textbf{223}, 509-532 (2001))
for Schr\"{o}dinger operators with sparse potential that increases as
fast as exponential--of--exponential. Applications to the spectrum of
the Kronecker sum of two (or more) copies of the model are given.
\end{abstract}

\section{Introduction}

\setcounter{equation}{0} \setcounter{theorem}{0}

The present paper deals with the Kronecker sum ($I$ denotes the identity
matrix) 
\begin{equation}
K=J\otimes I+I\otimes J  \label{K}
\end{equation}%
of off--diagonal Jacobi matrices%
\begin{equation}
J=\left( 
\begin{array}{ccccc}
0 & p_{0} & 0 & 0 & \cdots \\ 
p_{0} & 0 & p_{1} & 0 & \cdots \\ 
0 & p_{1} & 0 & p_{2} & \cdots \\ 
0 & 0 & p_{2} & 0 & \cdots \\ 
\vdots & \vdots & \vdots & \vdots & \ddots%
\end{array}%
\right) \;,  \label{J}
\end{equation}%
which are sparse perturbations of the free Jacobi matrix $J_{0}$ in the
sense that the sequence $(p_{n})_{n\geq 0}$ differs from the unit on a
lacunary subset of natural numbers $\mathbb{A}=\{a_{j}^{\omega }\}_{j\geq 1}$%
, i. e., 
\begin{equation}
p_{n}=\left\{ 
\begin{array}{ll}
p & \text{if~}n\in \mathbb{A}\,, \\ 
1 & \text{otherwise}\,,%
\end{array}%
\right.  \label{p}
\end{equation}%
\noindent for some $p\in (0,1)$. We assume that $\mathbb{A}$ is possibly a
random set, $a_{j}^{\omega }=a_{j}+\omega _{j}$, with%
\begin{equation*}
a_{j}^{\omega }-a_{j-1}^{\omega }\geq 2~
\end{equation*}%
and $a_{j}$ satisfying the "sparseness" condition $\lim_{j\rightarrow \infty
}a_{j}/a_{j-1}=\beta >1$ ($\beta =\infty $ is included). A concrete example
is given by 
\begin{equation}
a_{j}-a_{j-1}=\beta ^{j}\;,\qquad \qquad j=1,2,\ldots  \label{sparse}
\end{equation}%
\noindent with $a_{0}=0$, $\beta \geq 3$ and $\omega _{j}$, $j\geq 1$,
independent random variables defined on a probability space $(\Omega ,%
\mathcal{B},\mu )$, uniformly distributed over the sets $\Lambda _{j}\equiv
\left\{ -j,\ldots ,j\right\} $. These variables introduce uncertainty in the
position of the points where the $p_{n}$ differ from the unit and their
support increases linearly with the index $j$ (see Remark 1.4 of \cite{CMW1}
for less restrictive examples). A disordered potential of this type was
introduced by Zlato\v{s} \cite{Zla}. We shall however consider the
deterministic case with (\ref{sparse}) replaced by a sequence $\left( \beta
_{j}\right) _{j\geq 1}$ that increases faster--than--exponential.

The Jacobi matrices (\ref{J}) when applied to a vector $u=(u_{n})_{n\geq
0}\in l_{2}(\mathbb{Z}_{+})$ can be written as a difference equation 
\begin{equation*}
(Ju)_{n}=p_{n}u_{n+1}+p_{n-1}u_{n-1}\;,
\end{equation*}%
\noindent for $n\geq 0$ with $u_{-1}=0$. We denote by $J^{\phi }$ the Jacobi
matrix $J$ which satisfies $\phi $--boundary condition at $-1$: 
\begin{equation}
u_{-1}\cos \phi -u_{0}\sin \phi =0\;.  \label{phi}
\end{equation}%
$J^{\phi }$ is a (noncompact) perturbation of the free Jacobi matrix $%
J_{0}^{\phi }$, where $p_{n}=1$ for all $n\geq -1$: $J^{\phi }=J_{0}^{\phi
}+V$, the \textquotedblleft potential\textquotedblright\ $V$ composed by
infinitely many random barriers whose separations increase, at least,
exponentially fast. The Jacobi matrix $J$ corresponds to the matrix $J^{0}$
satisfying $0$--Dirichlet boundary condition at $-1$.

There have been a few results on models that exhibit spectral transition,
supporting spectra of different types in complementary set of parameters.
The Anderson model in a Bethe lattice \cite{K} is an example. Our intention
is to provide another instance of models whose spectrum contains pure point
and absolutely continuous nonempty components.  For the model (\ref{K}) with
random sparse potential satisfying (\ref{sparse}), the essential spectrum $%
\sigma _{\mathrm{ess}}(K)$ can be decomposed into continuous $\sigma _{%
\mathrm{c}}(K)$ at its center and dense pure point $\sigma _{\mathrm{pp}}(K)$
near to the edges. Whether $\sigma _{\mathrm{c}}(K)$ has an absolutely
continuous component still remains unknown.

A method of study the spectrum of sparse Jacobi matrices $J$, exploiting the
uniform distribution of the Pr\"{u}fer angles with fixed energy, has been
introduced in \cite{MarWre}. With this method, the Hausdorff dimension of
the spectral measure of sparse block--Jacobi matrices $J\otimes
I_{L}+I\otimes J_{L}$ can be determined with any degree of precision,
provided $\beta $ is large enough and $J_{L}$ homogeneous (see \cite{CMW}).
A sharp spectral transition from singular continuous to pure point spectrum,
announced by Zlato\v{s} for a model $J^{\phi }=J_{0}^{\phi }+V$ with
diagonal potential $V$ (see Theorem 6.3 of \cite{Zla}), has been proved in 
\cite{CMW1} for Jacobi matrices $J$ with random sparseness using this
method. Here we address the Fourier--Stieltjes transform $\widehat{d\rho }$
of the spectral measure $\rho $ of $J$ and investigate the least
faster--than--exponential sparse condition required for pointwise decay of $%
\widehat{d\rho }(t)$. While we consider off-diagonal case for reasons
explained in \cite{MarWre}, there will be no difficulties of principle in
applying our methods to the diagonal case considered in \cite{Zla}.

The layout of the present paper is as follows. Some notions and motivations
are presented in Section \ref{P} and precise statements, Theorem \ref%
{onehalf} and Corollary \ref{pac}, are formulated in Section \ref{S}.
Section \ref{S} also contains the mathematical tools used in Section \ref{PT}
to prove Theorem \ref{onehalf}. None of those tools depends on whether the Pr%
\"{u}fer angles are uniformly distributed. Lemma \ref{main} in Subsection %
\ref{BL} uses the idea of the proof of Theorem 10.12 in Chap. XII of \cite%
{Zy} and the Gevrey type estimates of Subsection \ref{GT} permit integration
by parts to be applied an unlimited number of times.

\section{Preliminaries\label{P}}

\setcounter{equation}{0} \setcounter{theorem}{0}

According to \cite{MarWre} (see e.g. Theorem 4.4), if the angles of Pr\"{u}%
fer are uniformly distributed, $\lambda \in \left[ -2,2\right] $ belongs to
the essential support of the singular continuous spectrum of $J$ provided%
\begin{equation}
\frac{r}{\beta }<1  \label{rbeta}
\end{equation}%
where $r=r(p,\lambda )=1+\vartheta (p)/(4-\lambda ^{2})$ and $\vartheta
(p)=\left( 1-p^{2}\right) ^{2}/p^{2}$ is monotone decreasing function of $%
p\in \left( 0,1\right) $. The local Hausdorff dimension of $\rho $ in this
case reads%
\begin{equation}
\alpha _{H}=\max \left( 1-\frac{\ln r}{\ln \beta },0\right) ~  \label{alpha}
\end{equation}%
(see Theorem 3.11 of \cite{CMW}). Note that $\alpha _{H}=\alpha _{H}(p,\beta
,\lambda )$ varies from $0$ to $1$ as $p$ varies from $p^{\ast }$ to $1$,
where $p^{\ast }$ is defined by $r(p^{\ast },\lambda )=\beta $; as $\lambda $
varies from $0$ to $\pm 2$, $\alpha _{H}$ varies from $\alpha _{H}(p,\beta
,0)<1$ to $0$, attained at $\lambda ^{\ast }$ defined by $r(p,\lambda ^{\ast
})=\beta $.

Let $f:\left[ -2,2\right] \longrightarrow \mathbb{C}$ be a smooth function
with $\text{supp}f\ni \lambda $ compact and sufficiently localized so that $%
\psi =f(J)\delta _{0}$ is in $l_{2}(\mathbb{Z}_{+})$. The probability $%
\left\vert \left( \psi ,\exp \left( itJ\right) \psi \right) \right\vert ^{2}$
of finding at time $t$ the system in its initial state $\psi $ can be
estimated observing that 
\begin{equation}
\left( \psi ,\exp \left( itJ\right) \psi \right) =\int_{-2}^{2}\left\vert
f(\lambda )\right\vert ^{2}e^{it\lambda }d\rho (\lambda )=\widehat{%
\left\vert f\right\vert ^{2}d\rho }(t)~,  \label{psi}
\end{equation}%
where $\rho (I)=\left\Vert E(I)\delta _{0}\right\Vert ^{2}$ is the spectral
measure of the state $\delta _{0}$ localized at $0$, and $E(\cdot )$ is the
spectral resolution of $J$. We shall also denote by $\rho _{\psi
}(I)=\left\Vert E(I)\psi \right\Vert ^{2}$ the spectral measure of an state $%
\psi =f(J)\delta _{0}$ and the Fourier--Stieltjes transform of $\rho _{\psi
} $ shall also be written as $\widehat{d\rho _{\psi }}(t)$.

Our purpose is primary to prove that (\ref{psi}) decays to zero as $%
\left\vert t\right\vert $ goes to infinity. Let $\alpha _{F}$ be the
supremum of $\alpha \geq 0$ such that 
\begin{equation}
\left\vert \widehat{|f|^{2}d\rho }(t)\right\vert \leq C\left( 1+\left\vert
t\right\vert \right) ^{-\alpha /2}~,  \label{alphaF}
\end{equation}%
with $C<\infty $, holds for any $t\in \mathbb{R}$, uniformly in a dense (in $%
L_{2}(-2,2)$) set of the $f$'s, sufficiently localized around a point $%
\lambda $ of the spectrum. It follows by a theorem of Frostman (see e.g. 
\cite{M}) that the so called Fourier dimension $\alpha _{F}$ satisfies%
\begin{equation*}
\alpha _{F}\leq \alpha _{H}
\end{equation*}%
(Fourier dimension and Hausdorff dimension do not agree in general).
Measures for which $\alpha _{F}>0$ and $\alpha _{F}=\alpha _{H}$ are,
respectively, known as Raychmann measures (see e.g. \cite{Ly}) and Salem
measures, after Salem's work \cite{S} on continuous distribution functions
which are constant in each interval contiguous to a perfect set of Lebesgue
measure zero.

Disorder plays a crucial role in most examples of sets and measures in which
the Hausdorff dimension and the Fourier dimension are equal (see e.g. \cite%
{M,Sa} and Chap. 17 of \cite{K}). The equality $\alpha _{H}=\alpha _{F}$ of
dimensions defined by (\ref{alpha}) and (\ref{alphaF}) is thus expected to
be attained for models with random sparseness condition (\ref{sparse}).
However, there are also singular continuous measures constructed by a
deterministic method (see e.g. \cite{Ko,P}) and for the deterministic sparse
model $J$ (or its diagonal version considered in references \cite{Zla,KR})
with faster--than--exponential ($\beta =\infty $) sparsity, it can actually
prove that this property is satisfied with $\alpha _{H}=\alpha _{F}=1$.

Pointwise decay 
\begin{equation}
\left\vert \widehat{\left\vert f\right\vert ^{2}d\rho }(t)\right\vert \leq
C\left\vert t\right\vert ^{-1/2}\ln \left\vert t\right\vert  \label{tlnt}
\end{equation}%
(with $\left\vert t\right\vert \geq 2$, because of the behavior of $\ln
\left\vert t\right\vert $ for $\left\vert t\right\vert \leq 1$) has been
obtained earlier by Simon \cite{S} for continuous Schr\"{o}dinger operators
with generic and sufficiently sparse potentials and, by using a different
method, Krutikov--Remling \cite{KR} have found, for a model similar to the
one considered here, a resonant set $\mathcal{R}\subset \mathbb{R}$ in which 
$\left\vert \widehat{\left\vert f\right\vert ^{2}d\rho }(t)\right\vert \leq
C\left( 1+\left\vert t\right\vert \right) ^{-1/2+\varepsilon }$ holds if $%
t\in \mathcal{R}$ and $\left\vert \widehat{\left\vert f\right\vert ^{2}d\rho 
}(t)\right\vert \leq C\left( 1+\left\vert t\right\vert \right) ^{-m}$ for an
arbitrary large but finite $m$, otherwise. The work of \cite{KR} was
motivated by the \textquotedblleft little control\textquotedblright\ of
reference \cite{S} on the rate with which the barrier separations $%
a_{j}-a_{j-1}$ have to increase (Simon's method requires $a_{j}\sim \exp
(\exp (cj^{3/2}))$ and Krutikov--Remling need e.g. $a_{j}\sim \exp (\exp
(cj))$, $c>0$, to satisfy their condition $a_{j}\leq Ca_{j+1}^{1-\mu }$ for
some $C>0$ and $\mu >0$).\footnote{%
The above mentioned decay $\left( 1+\left\vert t\right\vert \right)
^{-1/2+\varepsilon }$ holds in ref. \cite{KR} for every $t\in \mathcal{R}$
if $a_{j-1}<Ca_{j}^{1/2}$ is satisfied (i.e. for $\mu =1/2$) and $a_{j}\sim
\exp \left( c\exp j\right) $ with $c>\ln 2$. In ref. \cite{S}, this decay
holds provided $a_{j}\sim \exp \left( C_{\varepsilon }j^{3/2}\right) $ for
some $C_{\varepsilon }>0$.} In the present work, the sparseness condition is
less restrictive than those stated in references \cite{S} and \cite{KR}.
Thanks to the Gevrey type estimates developed in Subsection \ref{GT}, the
technique of integration by parts, employed in reference \cite{KR}, has been
exploited to its limit and (\ref{tlnt}) (with $\Omega (\left\vert
t\right\vert )$ in the place of $\ln \left\vert t\right\vert $, $\Omega (t)$
increasing less rapidly than any positive power of $t$) has been established
for $a_{j}-a_{j-1}$ increasing slightly--faster--than factorial: $\exp
\left( cj(\ln j)^{2}\right) /\delta ^{j}$ for some $c>1/2$ and $\delta <1$
(see (\ref{tepsilon}) for sparseness improvement: if $a_{j}-a_{j-1}=\exp
\left( \dfrac{1}{\varepsilon }j\ln j\right) /\delta ^{j}$ for some $%
\varepsilon >0$, then (\ref{tlnt}) holds with $\ln \left\vert t\right\vert $
replaced by $\left\vert t\right\vert ^{\varepsilon }$).

Our motivation in considering the Kronecker sum (\ref{K}) is the same that
led Simon to investigate the pointwise decay (\ref{tlnt}): it comes from the
observation that the convolution of two singular continuous measures, $d\rho
_{\psi }\ast d\rho _{\psi }$, may be absolutely continuous. Note that $%
\widehat{d\rho _{\psi }\ast d\rho _{\psi }}(t)=\left\vert \widehat{d\rho
_{\psi }}\right\vert ^{2}(t)$ is square integrable and this, according to
the well-known folklore result, implies that the corresponding measure is
absolutely continuous with respect to Lebesgue measure. For model with
faster--than--exponential sparseness the spectral measure is singular with
respect to Lebesgue and has Hausdorff dimension $1$, uniformly over the
essential spectrum $[-2,2]$ (see Theorem 1.4 of \cite{Zla}, which also
applies to the off--diagonal model). The spectrum of $K$ in this case is
purely absolutely continuous (see Corollary \ref{pac}). We are, however,
interested in spectral transitions and to achieve the decay (\ref{alphaF})
depending on the local Hausdorff dimension $\alpha _{H}$, a more
sophisticate method that exploits the randomness of $\{a_{j}^{\omega
}\}_{j\geq 1}$ is required. Investigation in this direction will be carried
in a separate paper.

As in \cite{KR}, our starting point is the representation of the spectral
measure as a weak--star--limit of absolutely continuous measures (see also 
\cite{P}):%
\begin{equation}
\int \mathrm{f}(\lambda )d\rho (\lambda )=\lim_{N\rightarrow \infty }\frac{1%
}{\pi }\int_{-2}^{2}\mathrm{f}(\lambda )\frac{\Im (w(\lambda ))}{\left\vert
y_{N}(\lambda )-w(\lambda )y_{N+1}(\lambda )\right\vert ^{2}}d\lambda
\label{int}
\end{equation}%
for every continuous function $\mathrm{f}:[-2,2]\longrightarrow \mathbb{C}$.
Here, $y_{n}=y_{n}(z)$ denotes the solution of the eigenvalue equation 
\begin{equation*}
p_{n}y_{n+1}+p_{n-1}y_{n-1}=zy_{n}~,
\end{equation*}%
satisfying the initial conditions $y_{-1}(z)=0$ and $y_{0}(z)=1$, and $w(z)=$
$z/2+i\sqrt{1-z^{2}/4}$ is a Herglotz function (maps the upper half-plane $%
\mathbb{H}$ into itself). We shall apply this formula with $\mathrm{f}%
(\lambda )=\left\vert f(\lambda )\right\vert ^{2}e^{it\lambda }$.

It turns out that $y_{n}(\lambda )$ is a subordinate solution of $Ju=\lambda
u$ (see details in \cite{CMW, CMW1}). Let $T(N,\lambda )$ denote the $%
2\times 2$ transfer matrix associated with the eigenvalue equation $%
Ju=\lambda u$ with $\lambda =2\cos \varphi $ for some $\varphi \in \left(
0,\pi \right) $ and define $\phi _{j}$, for the subsequence $\left(
N_{j}\right) _{j\geq 1}$ with $N_{j}=a_{j}+1$, by the equation 
\begin{equation*}
\left\vert T(N_{j},\lambda )\right\vert ^{2}\mathbf{v}_{\phi _{j}}=t_{j}^{-2}%
\mathbf{v}_{\phi _{j}}
\end{equation*}%
where $\left\vert T(N_{j},\lambda )\right\vert ^{2}=T^{\ast }(N_{j},\lambda
)T(N_{j},\lambda )$ is a real symmetric unimodular matrix, $t_{j}=\left\Vert
T(N_{j},\lambda )\right\Vert $ is the spectral norm of $T(N_{j},\lambda )$
and $\mathbf{v}_{\phi }=\dbinom{\cos \phi }{\sin \phi }$. Under the
assumption that the Pr\"{u}fer angles are uniformly distributed (satisfied
if $J$ is random; see \cite{CMW1}), $t_{j}^{2}=O(r^{j})$ and we have (see
Proposition 3.9 of \cite{CMW}) for an improved version of Lemma 2.1 of \cite%
{Zla}) 
\begin{equation}
\left\vert y_{N_{j}}(\lambda )-w(\lambda )y_{N_{j}+1}(\lambda )\right\vert
^{2}=\left\vert UT(N_{j},\lambda )\mathbf{v}_{\phi ^{\ast }}\right\vert
^{2}\ =O(r^{-j})\   \label{Rr}
\end{equation}%
for $\lambda $ in the essential support of $\rho $. Because the hypothesis
of Theorem 8.1 of Last--Simon \cite{LS} is verified, the limit $\phi ^{\ast
} $ of the sequence $\left( \phi _{j}\right) _{j\geq 1}$ exists and by the
Gilbert--Pearson theory $\phi ^{\ast }=0$ for a. e. $\varphi \in \left[
0,\pi \right] $ in the essential support of $\rho $ (see \cite{GP}),
establishing the decay property (\ref{Rr}).

We also have (with $\lambda =2\cos \varphi $ and $U=\left( 
\begin{array}{ll}
0 & \sin \varphi \\ 
1 & -\cos \varphi%
\end{array}%
\right) $) 
\begin{equation}
\left\vert UT(N_{j},\lambda )\mathbf{v}_{0}\right\vert ^{2}=R_{j}^{2}~,
\label{TR}
\end{equation}%
where $R_{j}=R_{j}(\varphi )$, $j=1,2,\ldots $, are the radius of Pr\"{u}fer
associated with $J$. Hence, our next ingredient is related with the
following identity%
\begin{equation}
\frac{1}{R_{j}^{2}}=\frac{1}{R_{j^{\ast }-1}^{2}}+\sum_{k=j^{\ast }-1}^{j-1}%
\frac{1}{R_{k}^{2}}\left( \frac{R_{k}^{2}}{R_{k+1}^{2}}-1\right)  \label{RR}
\end{equation}%
for some conveniently chosen $j^{\ast }\leq j$, together with%
\begin{equation}
\frac{R_{k}^{2}}{R_{k+1}^{2}}-1=\frac{p^{2}}{a+b\cos 2\theta _{k+1}+c\sin
2\theta _{k+1}}-1\equiv H(\varphi ,\theta _{k+1})  \label{RR-1}
\end{equation}%
where $\theta _{k}=\theta _{k}(\varphi )$ is the $k$--th angle of Pr\"{u}%
fer; $a$, $b$ and $c$ are functions of $p$ and $\varphi $ defined in \cite%
{MarWre} and $\bar{H}=\dfrac{1}{\pi }\displaystyle\int_{0}^{\pi }H(\varphi
,\theta )d\theta =0$. Plugging (\ref{RR}) and (\ref{RR-1}) into (\ref{int}),
yields

\begin{eqnarray}
\widehat{\left\vert f\right\vert ^{2}d\rho }(t) &=&\frac{1}{\pi }%
\int_{0}^{\pi }\left\vert f(2\cos \varphi )\right\vert ^{2}\frac{\sin
^{2}\varphi }{R_{j^{\ast }-1}^{2}}e^{2it\cos \varphi }d\varphi  \notag \\
&&+\sum_{k=j^{\ast }-1}^{\infty }\frac{1}{\pi }\int_{0}^{\pi }\left\vert
f(2\cos \varphi )\right\vert ^{2}\frac{\sin ^{2}\varphi }{R_{k}^{2}}%
H(\varphi ,\theta _{k+1})e^{2it\cos \varphi }d\varphi ~.~  \label{f2drho}
\end{eqnarray}

Expanding $H(\varphi ,\theta )$ in Fourier series (see \cite{KR}, for
details)%
\begin{eqnarray*}
H(\varphi ,\theta ) &=&\sum_{n=1}^{\infty }\left( A(\varphi
)^{n}e^{2in\theta }+\bar{A}(\varphi )^{n}e^{-2in\theta }\right) \\
A(\varphi ) &=&\sqrt{1-r^{-1}}e^{i(\delta +\pi )}~,
\end{eqnarray*}%
with $r$ given by (\ref{rbeta}) and $\tan \delta =c/b$,\footnote{%
Note that $\left\vert A(\varphi )\right\vert \leq a<1$ uniformly in each
compact set $K$ of $\left( 0,\pi \right) $ and the series is uniformly and
absolutely convergent.} each integral involved in the sum (\ref{f2drho}) is
of the form%
\begin{equation}
I(t)=\int_{-\infty }^{\infty }e^{ith(\varphi )}d\left( G\circ \lambda
\right) (\varphi )  \label{I}
\end{equation}%
for some $\mathcal{C}_{0}^{\infty }$ function $dG/d\lambda $ with support $%
2\cos [\varphi _{-},\varphi _{+}]\subset (-2,2)$, and $h$ of the type: 
\begin{equation*}
h_{k}(t,n;\varphi )=2(\cos \varphi +\frac{n}{t}\theta _{k}(\varphi ))\
,\qquad n\in \mathbb{Z}
\end{equation*}%
for some $k\geq j^{\ast }$.

The main contribution to (\ref{I}) comes from the stationary phase%
\begin{equation*}
h_{k}^{\prime }(\varphi )=-2\sin \varphi +\frac{2n}{t}\theta _{k}^{\prime
}(\varphi )=0
\end{equation*}%
and the $\varphi $'s that satisfy the equation will be called resonant or
critical values. Although the usual method of stationary phase does not
apply, its ideas can, nevertheless, be easily traced. To each integral for
which there are no resonant values inside the interval $\left[ \varphi
_{-},\varphi _{+}\right] $, the standard stationary phase method applies
integration by parts once. In order to apply integration by parts $m$ times,
with $m$ a large fixed number, obtaining therefore a pointwise decay $t^{-m}$%
, Krutikov--Remling in \cite{KR} have exploited the fact that $dG/d\lambda $
is analytic and the sequence $\left( a_{j}\right) _{j\geq 1}$ were
super-exponentially sparse.

We have in the present paper extended Krutikov--Remling's method in two
directions. Firstly, for pointwise decay $t^{-1}$ fixed to each of those no
resonant integrals, we use infinitely many integration by parts to obtain
the least possible sparseness condition. Secondly, concerning the integrals
in which there are resonant values inside the interval $\left[ \varphi
_{-},\varphi _{+}\right] $, we apply Lemma \ref{main}. Following the ideas
in the proof of Theorem 10.12 of Zygmund's book on trigonometric series, we
use Plancherel identity in order to obtain the main contribution of the
stationary phase, up to a logarithmic correction, provided the integral with
the phase $e^{ith(\varphi )}$ replaced by $e^{it\lambda (\varphi )}$ decays
as $t^{-1}$.

\section{Faster--than--exponential Sparse Models\label{S}}

\setcounter{equation}{0} \setcounter{theorem}{0}

\subsection{Statement of Results}

We devote this section to show that $\widehat{\left\vert f\right\vert
^{2}d\rho }(t)$ decays as $t^{-1/2}\Omega (t)$, uniformly with respect to
any $\mathcal{C}_{0}^{\infty }$ function $f$ with support contained into the
essential spectrum $[-2,2]$ of $J$, provided $\rho $ is the spectral measure
of $J$ given by (\ref{J}) and (\ref{p}) with the sparseness increment 
\begin{equation}
a_{j}-a_{j-1}=\beta _{j}\;,\qquad \qquad j=1,2,\ldots  \label{aa}
\end{equation}%
an increasing faster-than-exponential sequence: 
\begin{equation}
\lim_{j\rightarrow \infty }\beta _{j-1}/\beta _{j}=0.  \label{sparse1}
\end{equation}%
From now on, $J$ will always be such that (\ref{sparse1}) is satisfied. Our
goal is to find the least increasing sequence that leads to this result.
Without loss of generality, we assume that%
\begin{equation}
\frac{\beta _{j-1}}{\beta _{j}}\leq \delta  \label{delta}
\end{equation}%
holds uniformly in $j$ for some $0<\delta <1$ satisfying\footnote{%
By equation (4.15) of \cite{MarWre} $\delta >\min_{\varphi ^{-}\leq \varphi
\leq \varphi ^{+}}\left( a/p-\sqrt{(a/p)^{2}-1}\right) >0$ if $0<p<1$.}%
\begin{equation}
\sup_{\varphi \in \text{supp}f\circ \lambda }\sup_{\theta }\frac{p^{2}}{%
a+b\cos 2\theta +c\sin 2\theta }<\frac{1}{\delta }~  \label{delta-l}
\end{equation}%
Note that%
\begin{equation*}
\frac{\delta ^{j}}{R_{j}^{2}(\varphi )}\longrightarrow 0
\end{equation*}%
exponentially fast, as $j$ goes to infinity, uniformly in supp$f\circ
\lambda $. We now state our result.

\begin{theorem}
\label{onehalf}Let $\rho $ the spectral measure of $J$ associated with the
state $\delta _{0}$ localized at $0$ and let $f$ be a smooth function with
compact support inside $\left( -2,2\right) $ and such that $0\notin \text{%
supp}f$. Suppose that the sequence $\left( a_{j}\right) _{j\geq 1}$
satisfies (\ref{aa})--(\ref{delta-l}) with 
\begin{equation}
\beta _{j}=\frac{1}{\delta ^{j}}\exp \left( cj\left( \ln j\right) ^{2}\right)
\label{betaj}
\end{equation}%
for some $c>1/2$, as $j$ tends to infinity. Then, there exist a constant $C$%
, depending on $f$ and $p$ (the intensity parameter of $J$), such that%
\begin{equation}
\left\vert \widehat{\left\vert f\right\vert ^{2}d\rho }(t)\right\vert \leq
C\left\vert t\right\vert ^{-1/2}~\Omega (\left\vert t\right\vert )~
\label{rhot1/2}
\end{equation}%
holds for $\left\vert t\right\vert \geq 2$, where $\Omega (t)$ increases
less rapidly than any positive power of $t$.

Moreover, if the sparseness increments $\left( \beta _{j}\right) _{j\geq 1}$
is chosen as 
\begin{equation}
\beta _{j}=\frac{1}{\delta ^{j}}\exp \left( \varepsilon ^{-1}j\ln j\right)
\label{factorial}
\end{equation}%
for some $\varepsilon >0$ small enough, then the conclusion (\ref{rhot1/2})
holds with the upper bound replaced by $C(1+\left\vert t\right\vert
)^{-1/2+\varepsilon }$.
\end{theorem}

The proof of Theorem \ref{onehalf} uses a classical result on the decay of
Fourier--Stieltjes coefficients $c_{n}(dG)$ of a monotone increasing
singular continuous function $G$ originated from a Riesz product (see
Theorem $10.12$ in Chap. XII of \cite{Zy}). The Lemma stated below reduces
the resonant estimate to a non--resonant one, which will be studied in
Subsection \ref{PT}. Our proof resembles, in this sense, the proof of Simon
(compare Lemma \ref{main} below with Lemma $4.2$ of \cite{S}) with the
non--resonant estimate given by Krutikov--Remling's method of integration by
parts.

A well-known folklore result (see, e.g., \cite{C}, Exercise 11, Section 6.2)
states that if the Fourier-Stieltjes transform of a finite Borel measure is
square-integrable, the corresponding measure is absolutely continuous with
respect to Lebesgue measure. We have, thus, as a direct corollary of Theorem %
\ref{onehalf}:

\begin{corollary}
\label{pac}The spectrum of $K=J\otimes I+I\otimes J$, with $J$ defined as in
Theorem \ref{onehalf}, is purely absolutely continuous.
\end{corollary}

An alternative proof of a theorem that includes the above mentioned folklore
result (which has been used as early as 1958 to provide nontrivial examples
of the statement that the convolution of two singular measures may be
absolutely continuous; see \cite{KS}) may be found in \cite{CMW2} (see also 
\cite{S}, Corollary 3.2). Note that the spectral measure of $K$ associated
with the tensor product of two vectors in $l_{2}(\mathbb{Z}_{+})$ is given
by the convolution of the two measures associated to each of these vectors
(see \cite{S} for details).

Of course the corollary extends to the Kronecker sum of any number of copies
higher than two, and thus the nature of the spectrum changes dramatically in
dimension two or higher. This will be exploited further in \cite{CMW2}.

\subsection{Basic Lemma\label{BL}}

The goal of this subsection is to prove the following

\begin{lemma}
\label{main}Suppose $G:\mathbb{R}\longrightarrow \mathbb{R}$ is a monotone
increasing continuous function, with $dG$ supported in some closed interval $%
[a,b]\subset \left( -2,0\right) \cup (0,2)$, whose Fourier--Stieltjes
transform satisfies%
\begin{equation}
\widehat{dG}(t)=\int_{-\infty }^{\infty }e^{it\lambda }dG(\lambda )\leq 
\frac{C}{1+\left\vert t\right\vert }  \label{Ghat}
\end{equation}%
for some constant $C<\infty $ and every $t\in \mathbb{R}$. Let $\gamma :%
\mathbb{R}\longrightarrow \mathbb{C}$ be defined by%
\begin{equation*}
\gamma (t)=\int_{-\infty }^{\infty }e^{itx(t,\lambda )}dG(\lambda )
\end{equation*}%
where%
\begin{equation}
tx(t,\lambda )=t\lambda +\frac{\kappa }{\pi }\cos ^{-1}\frac{\lambda }{2}
\label{tx}
\end{equation}%
is a mapping from $[-2t,2t]$ into $[-2t+\kappa ,2t]$. If $\kappa =\kappa
(t)=O\left( \left\vert t\right\vert \right) $, then 
\begin{equation*}
\left\vert \gamma (t)\right\vert \leq \frac{B}{\left\vert t\right\vert ^{1/2}%
}\ln \left\vert t\right\vert
\end{equation*}%
holds for some $B<\infty $ and every $t\in \mathbb{R}$ with $\left\vert
t\right\vert \geq 2$.
\end{lemma}

\noindent \textit{Proof.} Denoting by $\chi $ the characteristic function of
the interval $\left[ a,b\right] $: $\chi (\lambda )=1$ if $a\leq \lambda
\leq b$ and $=0$ otherwise, by the Plancherel theorem%
\begin{equation}
\gamma (t)=\int_{-\infty }^{\infty }e^{itx(t,\lambda )}\chi (\lambda
)dG(\lambda )=\int_{-\infty }^{\infty }\Lambda (t,\tau )\widehat{dG}(\tau
)d\tau  \label{integral}
\end{equation}%
where%
\begin{equation}
\Lambda (t,\tau )=\frac{1}{2\pi }\int_{a}^{b}e^{i\left( tx(t,\lambda )+\tau
\lambda \right) }d\lambda ~.  \label{lambdahat}
\end{equation}%
We apply van der Corput estimates (see e.g. Lemma $4.3$, Chap. $V$ of \cite%
{Zy}) in order to obtain the asymptotic behavior of $\Lambda (t,\tau )$ for
large $t$ and $\tau $. The integral (\ref{lambdahat}) is of the form%
\begin{eqnarray*}
2\pi \Lambda (t,\tau ) &=&\int_{a}^{b}e^{2\pi if(\lambda )}d\lambda \\
2\pi f(\lambda ) &=&(t+\tau )\lambda +\frac{\kappa }{\pi }\cos ^{-1}\frac{%
\lambda }{2}
\end{eqnarray*}%
where the second derivative $f^{\prime \prime }$ of $f$ is strictly negative
(under the hypothesis $0\notin \lbrack a,b]$, $\cos ^{-1}\lambda /2$ is
strictly concave for every $\lambda \in \left[ a,b\right] $) and
proportional to $\kappa =O\left( t\right) $, uniformly in $\tau $: $%
\left\vert f^{\prime \prime }(\lambda )\right\vert \geq \rho \left\vert
\kappa \right\vert $ where $2\pi ^{2}\rho =\sup_{\lambda \in \left[ a,b%
\right] }\lambda /(4-\lambda ^{2})^{3/2}$. By van der Corput's lemma%
\begin{equation}
\left\vert \Lambda (t,\tau )\right\vert \leq \frac{4}{\sqrt{\rho \kappa }}%
\leq \frac{K}{\sqrt{1+\left\vert t\right\vert }}  \label{K1}
\end{equation}%
holds for a constant $K$, independent of $t$ and $\tau $. On the other hand,
for $\left\vert \tau \right\vert $ large compared with $\left\vert
t\right\vert $, let us say $\left\vert \tau \right\vert >\Delta \left\vert
t\right\vert $ for $\Delta =1+\left\vert \kappa /t\right\vert \sup_{\lambda
\in \lbrack a,b]}1\left/ \left( \pi \sqrt{4-\lambda ^{2}}\right) \right. $,
the derivative $f^{\prime }$ of $f$ is of order $\tau $ and Van der Corput's
lemma gives%
\begin{equation}
\left\vert \Lambda (t,\tau )\right\vert \leq \frac{K^{\prime }}{\left\vert
\tau \right\vert }  \label{K2}
\end{equation}%
for some constant $K^{\prime }$, independent of $t$ and $\tau $.

Estimates (\ref{K1}) and (\ref{K2}) together with (\ref{Ghat}), yield 
\begin{eqnarray}
\left\vert \gamma (t)\right\vert &\leq &\int_{\left\vert \tau \right\vert
\leq \Delta \left\vert t\right\vert }\left\vert \Lambda (t,\tau )\right\vert
\left\vert \widehat{dG}(\tau )\right\vert d\tau +\int_{\left\vert \tau
\right\vert >\Delta \left\vert t\right\vert }\left\vert \Lambda (t,\tau
)\right\vert \left\vert \widehat{dG}(\tau )\right\vert d\tau  \notag \\
&\leq &\frac{K}{\sqrt{1+\left\vert t\right\vert }}\int_{\left\vert \tau
\right\vert \leq \Delta \left\vert t\right\vert }\frac{C}{1+\left\vert \tau
\right\vert }d\tau +\int_{\left\vert \tau \right\vert >\Delta \left\vert
t\right\vert }\frac{K^{\prime }}{\left\vert \tau \right\vert }\frac{C}{%
1+\left\vert \tau \right\vert }d\tau  \notag \\
&\leq &B^{\prime }\left\vert t\right\vert ^{-1/2}\ln \left\vert t\right\vert
+B^{\prime \prime }\left\vert t\right\vert ^{-1}  \label{dF}
\end{eqnarray}%
for $\left\vert t\right\vert \geq 2$ and some finite constants $B^{\prime }$
and $B^{\prime \prime }$, concluding the proof of the lemma.

\hfill $\Box $

\begin{remark}
\label{tepsilon}If the pointwise behavior of (\ref{Ghat}) is replaced by $%
C/(1+\left\vert t\right\vert )^{1-\varepsilon }$ for some $\varepsilon >0$,
then the logarithmic correction $\ln \left\vert t\right\vert $ in the
conclusion of Lemma has to be replaced by $\left\vert t\right\vert
^{\varepsilon }$ (see the last inequality of (\ref{dF})). We shall also
apply Lemma \ref{main} with $n^{\ast }\kappa =O(t)$ where $n^{\ast }=n^{\ast
}(t)$ is a nonvanishing integer valued, piecewise monotone function,
increasing slower than any positive power of $t$. In this case, logarithmic
correction will be replaced by $\sqrt{n^{\ast }(t)}$ (see Section \ref{PT},
for details).
\end{remark}

\begin{remark}
The function $x=x(\lambda )$ in Theorem 10.12 of Chap. XII of \cite{Zy},
whose proof has suggested us Lemma \ref{main}, is a one--to--one mapping of $%
\left[ -\pi ,\pi \right] $ onto itself (it does not depend on $t$, as in our
case; see (\ref{tx})) and $\gamma (t)=\widehat{dF}(t)$, thereby, where $F(x)$
is an increasing function with $F\circ x(\lambda )=\sqrt{2\pi }G(\lambda )$.
These facts are not necessary for the conclusion of Lemma \ref{main}. The
hypothesis on $\widehat{dG}(t)$ in that theorem is, in addition, stronger
than ours. Because $\widehat{dG}(t)$ decays only \textquotedblleft on the
average\textquotedblright\ it is necessary an improved estimate $K^{\prime
}/\tau ^{2}$ instead of a simpler one (\ref{K2}).
\end{remark}

\subsection{Gevrey Type Estimates\label{GT}}

To prove Theorem \ref{onehalf} we need to improve results of \cite{KR} in
order to bring $\beta _{j}$ down from the asymptote $\exp \left( \exp
cj\right) $ to $\left( j!\right) ^{1/\varepsilon }$ for some $\varepsilon >0$%
. The following proposition gives a Gevrey type estimate for the derivatives
of the Pr\"{u}fer angles.

\begin{proposition}
\label{scott}Let $\theta _{k}=\theta _{k}(\varphi )$, $k=0,1,\ldots $, be
the sequence of Pr\"{u}fer angles starting from $\theta _{0}$, consistent
with the initial condition (\ref{Rr}). Then, for every $m\geq 1$, 
\begin{equation}
\left\vert \theta _{m}^{\prime }(\varphi )\right\vert \leq C_{1}\delta \eta
\beta _{m}  \label{eta}
\end{equation}%
and for every $m>1$ and $n>1$ 
\begin{equation}
\frac{1}{n!}\left\vert \theta _{m}^{(n)}(\varphi )\right\vert \leq C_{n}\eta
^{n}\beta _{m-1}^{n}~  \label{ineq}
\end{equation}%
hold uniformly in compact subsets of $\left( 0,\pi \right) $ with $%
C_{n}=K/n^{2}$, $K\leq $ $3/(2\pi ^{2})=0.151981...$ and 
\begin{equation*}
\eta =\frac{1+\Delta }{\delta K}
\end{equation*}%
where $\Delta <1$ ($\Delta =O(\delta )$) is a constant satisfying $\theta
_{m}^{\prime }(\varphi )\leq \left( 1+\Delta \right) \beta _{m}$, which is
computable from Proposition 5.2 of \cite{MarWre}, and $\delta $ is suitably
small.
\end{proposition}

\begin{remark}
Proposition \ref{scott} replaces the unspecified constant $C_{j}$ appearing
in Lemma $3.1$ of \cite{KR} by $K\eta ^{j}j!/j^{2}$. Detailed information on
the growth of $\theta _{m}^{(j)}(\varphi )$ in both $j$ and $m$ are an
essential ingredient of our method.
\end{remark}

To prove Proposition \ref{scott} we need the following

\begin{lemma}
\label{A}Let $\mathbf{C}\ast \mathbf{D}$ denote the convolution product in $%
\mathbb{R}^{\mathbb{Z}_{+}}$:%
\begin{equation}
\left( \mathbf{C}\ast \mathbf{D}\right) _{n}=\sum_{i=0}^{n}C_{i}D_{n-i}
\label{prod}
\end{equation}%
for $n\geq 0$ . If $\mathbf{C}$ has components given by $C_{i}=K/i^{2}$, $%
i\geq 1$, with $C_{0}=K$, then%
\begin{equation*}
\underset{k-\text{factors}}{\underbrace{\mathbf{C}\ast \mathbf{C}\ast \cdots
\ast \mathbf{C}}}~\leq \mathbf{C}
\end{equation*}%
holds for every $k\geq 1$ provided $K\leq 1/\left( 2+2\pi ^{2}/3\right) $.
If $C_{0}=0$, then same result holds with $K\leq 3/(2\pi ^{2})$.
\end{lemma}

\begin{remark}
The $n$--th component of the convolution of $\mathbf{C}$ with itself $k$%
--times satisfies%
\begin{equation}
\sum_{\substack{ i_{1},\ldots ,i_{k}\geq 0  \\ i_{1}+\cdots +i_{k}=n}}%
C_{i_{1}}\cdots C_{i_{k}}\leq C_{n}  \label{C}
\end{equation}%
by Lemma \ref{A}.
\end{remark}

\begin{remark}
\label{variation}We have stated Lemma \ref{A} for sequences $\mathbf{C}%
=\left( C_{0},C_{1},\ldots \right) \in \mathbb{R}^{\mathbb{Z}_{+}}$ with the 
$0$--th component $C_{0}=0$ and $C_{0}=K$ since both cases will be
considered in this subsection (see Proposition \ref{estimates} below for the
case $C_{0}\neq 0$). Lemma \ref{A} plays a key role in every estimate
involving higher order chain rule.
\end{remark}

\noindent \textit{Proof of lemma.} By definition,%
\begin{equation}
\frac{1}{C_{n}}\sum_{i=0}^{n}C_{i}C_{n-i}=K\left( 2+\sum_{i=1}^{n-1}\frac{%
n^{2}}{i^{2}(n-i)^{2}}\right) ~.  \label{CCC}
\end{equation}%
Writing

\begin{equation}
\frac{n^{2}}{i^{2}(n-i)^{2}}=\frac{n^{2}}{\left( n-i\right) ^{2}+i^{2}}%
\left( \frac{1}{i^{2}}+\frac{1}{\left( n-i\right) ^{2}}\right) ~,  \label{1}
\end{equation}%
the pre--factor in the r.h.s. of (\ref{1}) can be bounded using $0\leq
\left( a-b\right) ^{2}=2\left( a^{2}+b^{2}\right) -\left( a+b\right) ^{2}$,
which holds for any real numbers $a$ and $b$, with $a=n-i$ and $b=i$. We have%
\begin{equation*}
0\leq (n-2i)^{2}=2\left( (n-i)^{2}+i^{2}\right) -n^{2}
\end{equation*}%
or, equivalently,%
\begin{equation}
\frac{n^{2}}{\left( n-i\right) ^{2}+i^{2}}\leq 2~.  \label{2}
\end{equation}%
Plugging (\ref{1}) into (\ref{CCC}) together with (\ref{2}), gives 
\begin{eqnarray}
\frac{1}{C_{n}}\sum_{i=0}^{n}C_{i}C_{n-i} &\leq &2K\left(
1+\sum_{i=1}^{n-1}\left( \frac{1}{i^{2}}+\frac{1}{(n-i)^{2}}\right) \right) 
\notag \\
&\leq &2K\left( 1+\frac{\pi ^{2}}{3}\right) \leq 1  \label{CC/C}
\end{eqnarray}%
provided $K\leq 1/\left( 2+2\pi ^{2}/3\right) $. The case $C_{0}=0$, the
terms with $i=0$ and $n$ do not contribute to the sum and the inequality (%
\ref{CC/C}) holds provided $K\leq 3/(2\pi ^{2})$. Once we have $\mathbf{C}%
\ast \mathbf{C}~\leq \mathbf{C}$, Lemma \ref{A} is proved by induction.

\hfill $\Box $

\noindent \textit{Proof of Proposition \ref{scott}.} The proof uses the
recursive relation%
\begin{equation}
\theta _{m}(\varphi )=g\circ \theta _{m-1}(\varphi )+\beta _{m}\varphi
\label{thetatheta}
\end{equation}%
where 
\begin{equation}
g=g(\varphi ,\theta )=\tan ^{-1}\left( (\tan \theta +\cot \varphi
)/p^{2}-\cot \varphi \right)  \label{gg}
\end{equation}%
together with the Scott's formula for higher order chain rule (see e.g. \cite%
{FLy})%
\begin{equation}
\left( g\circ f\right) ^{[n]}=\sum_{k=1}^{n}g^{[k]}\circ f\sum_{\substack{ %
i_{1},\ldots ,i_{k}\geq 1  \\ i_{1}+\cdots +i_{k}=n}}f^{[i_{1}]}\cdots
f^{[i_{k}]}  \label{sf}
\end{equation}%
where, from here on, $h^{[n]}$ stands for $h^{(n)}/n!$, the $n$--th
derivative of $h$ divided by $n!$.

Upper and lower bounds for the first derivative has been provided in \cite%
{MarWre}:%
\begin{equation}
\left( 1-\Delta \right) \beta _{m}\leq \theta _{m}^{\prime }(\varphi )\leq
\left( 1+\Delta \right) \beta _{m}  \label{beta-beta}
\end{equation}%
with $\Delta <1$ a constant. Now, choosing $\eta =\dfrac{1+\Delta }{\delta K}
$, (\ref{beta-beta}) establishes (\ref{eta}) for every $m\geq 1$.

Since $\left( \beta _{j}\right) _{j\geq 1}$ is a fast increasing sequence we
apply the Scott's formula to $g$ in (\ref{gg}) as it were a function of a
single variable $\theta $. This really gives the main contribution to the
derivatives. $g$ as a function of $z=e^{i\theta }$, continued to the complex
plane, is analytic outside a disc of radius strictly less than $1-e/\xi <1$.
The derivatives of $q(e^{i\theta })=g(\theta )$ may be estimate by Cauchy
formula:%
\begin{equation}
\left\vert g^{[k]}(\theta )\right\vert \leq c_{1}\xi ^{k}~  \label{g}
\end{equation}%
holds for $k\geq 1$ with $c_{1}$ as small as one wishes, by increasing $\xi $
accordingly ((\ref{g}) can be bounded, e. g., by $\varepsilon \left(
c_{1}\xi /\varepsilon \right) ^{k}=\varepsilon \bar{\xi}^{k}$, for any $%
\varepsilon >0$). Replacing $f$ by $\theta _{m-1}$ in (\ref{sf}), gives 
\begin{equation}
\theta _{m}^{[n]}(\varphi )=\sum_{k=1}^{n}g^{[k]}\circ \theta _{m-1}(\varphi
)\sum_{\substack{ i_{1},\ldots ,i_{k}\geq 1  \\ i_{1}+\cdots +i_{k}=n}}%
\theta _{m-1}^{[i_{1}]}(\varphi )\cdots \theta _{m-1}^{[i_{k}]}(\varphi )~.
\label{theta}
\end{equation}

We prove (\ref{ineq}) by induction in $n$. Consider the case $n=2$, for any $%
m>1$. By equation (\ref{thetatheta}), together with (\ref{g}) and (\ref{eta}%
), we have%
\begin{eqnarray*}
\theta _{m}^{\prime \prime } &=&g^{\prime \prime }\circ \theta _{m-1}\cdot
(\theta _{m-1}^{\prime })^{2}+g^{\prime }\circ \theta _{m-1}\cdot \theta
_{m-1}^{\prime \prime } \\
&\leq &c_{1}\xi ^{2}C_{1}^{2}\delta ^{2}\eta ^{2}\beta _{m-1}^{2}+g^{\prime
}\circ \theta _{m-1}\cdot \theta _{m-1}^{\prime \prime }~.
\end{eqnarray*}%
The iteration of this relation together with (\ref{g}) and (\ref{delta}),
yields%
\begin{eqnarray*}
\theta _{m}^{\prime \prime } &\leq &c_{1}\xi ^{2}C_{1}^{2}\delta ^{2}\eta
^{2}\beta _{m-1}^{2}\sum_{j=1}^{m-1}\left( c_{1}\xi \right) ^{j-1}\frac{%
\beta _{m-j}^{2}}{\beta _{m-1}^{2}} \\
&\leq &\frac{c_{1}\xi ^{2}\delta ^{2}}{1-c_{1}\xi \delta ^{2}}C_{1}^{2}\eta
^{2}\beta _{m-1}^{2} \\
&\leq &2C_{2}\eta ^{2}\beta _{m-1}^{2}
\end{eqnarray*}%
provided $\delta $ is chosen so small that $c_{1}\xi \delta ^{2}<1$ and%
\begin{equation*}
2K\frac{c_{1}\xi ^{2}\delta ^{2}}{1-c_{1}\xi \delta ^{2}}\leq 1~
\end{equation*}%
are both satisfied, establishing (\ref{ineq}) for $n=2$.

Now, suppose%
\begin{equation*}
\theta _{m}^{[j]}(\varphi )\leq C_{j}\eta ^{j}\beta _{m-1}^{j}
\end{equation*}%
holds for $m>1$ and $j=2,\ldots ,n-1$ and we shall establish the inequality
for $n$. By this assumption together with (\ref{delta}), we have%
\begin{equation}
\theta _{m-1}^{[j]}(\varphi )\leq C_{j}\eta ^{j}\beta _{m-2}^{j}=C_{j}\eta
^{j}\left( \frac{\beta _{m-2}}{\beta _{m-1}}\right) ^{j}\beta _{m-1}^{j}\leq
C_{j}(\delta \eta )^{j}\beta _{m-1}^{j}~.  \label{delta-eta}
\end{equation}%
Plugging (\ref{eta}), (\ref{g}) and (\ref{delta-eta}) into (\ref{theta}),
together with (\ref{C}), yields%
\begin{eqnarray*}
\theta _{m}^{[n]}(\varphi ) &\leq &c_{1}(\delta \eta )^{n}\beta
_{m-1}^{n}\sum_{k=2}^{n}\xi ^{k}\sum_{\substack{ i_{1},\ldots ,i_{k}\geq 1 
\\ i_{1}+\cdots +i_{k}=n}}C_{i_{1}}\cdots C_{i_{k}}+g^{\prime }\circ \theta
_{m-1}\cdot \theta _{m-1}^{[n]} \\
&\leq &c_{1}\frac{\xi }{\xi -1}C_{n}(\delta \xi \eta )^{n}\beta
_{m-1}^{n}+g^{\prime }\circ \theta _{m-1}\cdot \theta _{m-1}^{[n]}
\end{eqnarray*}%
Here, we have separated the term with $k=1$ which applies $n$ derivatives on 
$\theta _{m-1}$. Note that, for all the other terms with $k\geq 2$, we have $%
i_{1},\ldots ,i_{k}\geq 1$ and the derivatives applied on the $\theta _{m-1}$
are of order strictly smaller than $n$. The iteration of this relation,
gives 
\begin{eqnarray*}
\theta _{m}^{[n]}(\varphi ) &\leq &c_{1}\frac{\xi }{\xi -1}C_{n}(\delta \xi
\eta )^{n}\beta _{m-1}^{n}\sum_{j=1}^{m-1}\left( c_{1}\xi \right) ^{j-1}%
\frac{\beta _{m-j}^{n}}{\beta _{m-1}^{n}} \\
&\leq &C_{n}\eta ^{n}\beta _{m-1}^{n}
\end{eqnarray*}%
provided%
\begin{equation}
c_{1}\frac{\xi }{\xi -1}\left( \xi \delta \right) ^{n}\frac{1}{1-c_{1}\xi
\delta ^{n}}\leq 1  \label{c1}
\end{equation}%
holds for every $n>2$. We pick $\delta $ satisfying both (\ref{delta}) and (%
\ref{c1}) for $n\geq 2$, concluding the proof of Proposition \ref{scott}.

\hfill $\Box $

\begin{remark}
The well known formula for higher derivative of composite functions, Fa\`{a}
di Bruno's formula, cannot be used recursively since the constant $\eta $ in
equation (\ref{ineq}) deteriorates each time it is applied (see eq. (6.10)
in Subsection 6.2 of \cite{BM}). The proof of (\ref{ineq}), by induction,
using Scott's formula was based on yet unpublished manuscript
\textquotedblleft $O(N)$ Hierarchical Model Approached by the Implicit
Function Theorem\textquotedblright\ by W. R. P. Conti and D. H. U. Marchetti.
\end{remark}

To obtain the $t^{-1}$ decay from the summation in (\ref{f2drho}), it is
necessary to apply an arbitrarily large number of the integration by parts
for integrals of the type (\ref{I}):%
\begin{equation*}
\int_{0}^{\pi }f_{0}e^{ith}d\varphi =i\int_{0}^{\pi }\left( \frac{1}{%
th^{\prime }}f_{0}\right) ^{\prime }e^{ith}d\varphi
\end{equation*}%
where $\text{supp}f_{0}=[\varphi _{-},\varphi _{+}]\subset \left( 0,\pi
\right) $. The following propositions gather tools to implement the estimate.

\begin{proposition}
\label{Ln}Let $f_{0}(\varphi )$ and $\varrho (\varphi )$ be, respectively, $%
\mathcal{C}^{\infty }$ complex and real--valued functions on $[0,\pi )$ and
let $L=\dfrac{d}{d\varphi }\varrho (\varphi )$ be an operator defined in
this space. If%
\begin{equation*}
f_{n}=\frac{1}{n!}L^{n}f_{0}=\frac{1}{n!}\dfrac{d}{d\varphi }\varrho \dfrac{d%
}{d\varphi }\varrho \cdots \dfrac{d}{d\varphi }\varrho f_{0}
\end{equation*}%
denotes the $n$--th application of $L$ over $f_{0}$, divided by factorial of 
$n$, $n=0,1,\ldots $, then%
\begin{equation}
f_{n}=\sum_{_{\substack{ k_{1},\ldots ,k_{n},p_{n}\geq 0  \\ k_{1}+\cdots
+k_{n}+p_{n}=n}}}\varrho ^{\lbrack k_{1}]}\cdots \varrho ^{\lbrack
k_{n}]}f_{0}^{[p_{n}]}~.  \label{f0}
\end{equation}
\end{proposition}

\noindent \textit{Proof.} The proof is by induction. For $n=1$, we have $%
f_{1}=\varrho ^{\prime }f_{0}+\varrho f_{0}^{\prime }$. Assuming that (\ref%
{f0}) holds, 
\begin{eqnarray*}
f_{n+1} &=&\frac{1}{n+1}\left( \varrho ^{\prime }f_{n}+\varrho f_{n}^{\prime
}\right) \\
&=&\frac{1}{n+1}\sum_{\substack{ \hat{k}_{1},\ldots ,\hat{k}_{n+1},\hat{p}%
_{n+1}\geq 0  \\ \hat{k}_{1}+\cdots +\hat{k}_{n+1}+\hat{p}_{n+1}=n+1}}(\hat{k%
}_{1}+\cdots +\hat{k}_{n+1}+\hat{p}_{n+1})\varrho ^{\lbrack \hat{k}%
_{1}]}\cdots \varrho ^{\lbrack \hat{k}_{n+1}]}f_{0}^{[\hat{p}_{n+1}]} \\
&=&\sum_{\substack{ \hat{k}_{1},\ldots ,\hat{k}_{n+1},\hat{p}_{n+1}\geq 0 
\\ \hat{k}_{1}+\cdots +\hat{k}_{n+1}+\hat{p}_{n+1}=n+1}}\varrho ^{\lbrack 
\hat{k}_{1}]}\cdots \varrho ^{\lbrack \hat{k}_{n+1}]}f_{0}^{[\hat{p}_{n+1}]}
\end{eqnarray*}%
where we have applied the product rule%
\begin{equation*}
f_{n}^{\prime }=\sum_{j=1}^{n}\sum_{_{\substack{ k_{1},\ldots
,k_{n},p_{n}\geq 0  \\ k_{1}+\cdots +k_{n}+p_{n}=n}}}\varrho ^{\lbrack
k_{1}]}\cdots \left( \varrho ^{\lbrack k_{j}]}\right) ^{\prime }\cdots
\varrho ^{\lbrack k_{n}]}f_{0}^{[p_{n}]}+\sum_{_{\substack{ k_{1},\ldots
,k_{n},p_{n}\geq 0  \\ k_{1}+\cdots +k_{n}+p_{n}=n}}}\varrho ^{\lbrack
k_{1}]}\cdots \varrho ^{\lbrack k_{n}]}\left( f_{0}^{[p_{n}]}\right)
^{\prime }~,
\end{equation*}%
used $\left( \varrho ^{\lbrack k_{j}]}\right) ^{\prime }=(k_{j}+1)\varrho
^{\lbrack k_{j}+1]}$ for $j=1,\ldots ,n$ (analogously for $\left(
f_{0}^{[p_{n}]}\right) ^{\prime }=(p_{n}+1)\varrho ^{\lbrack p_{n}+1]}$),
redefined variables: 
\begin{equation*}
\hat{k}_{j+1}=k_{j}+1~,\qquad \hat{k}_{l+1}=k_{l}\ \ \text{for }l\neq
j\qquad \text{and}\qquad \hat{p}_{n+1}=p_{n}
\end{equation*}%
(the same for $\hat{p}_{n+1}=p_{n}+1$ and $\hat{k}_{j+1}=k_{j}$ for $%
j=1,,\ldots ,n$) and have added a new variable $\hat{k}_{1}$. Observe that $%
\hat{k}_{j+1}=k_{j}+1\geq 1$ but we can start the sum over $\hat{k}_{j+1}$
from $0$ since $\hat{k}_{j+1}\varrho ^{\lbrack \hat{k}_{j+1}]}$ is
identically $0$ at $\hat{k}_{j+1}=0$. This completes the proof of the
proposition.

\hfill $\Box $

In our application, $\varrho =\dfrac{1}{th_{m+1}^{\prime }(\varphi )}$ and $%
f_{0}=\left\vert f(2\cos \varphi )\right\vert ^{2}\dfrac{\sin ^{2}\varphi }{%
R_{m}^{2}(\varphi )}A^{n}(\varphi )$ (or its complex conjugate). The main
contribution for $m\geq j^{\ast }+1$, where $j^{\ast }=j^{\ast }(t)$ is such
that%
\begin{equation}
\beta _{j^{\ast }}\leq t<\beta _{j^{\ast }+1}~,  \label{t}
\end{equation}%
comes from the derivatives of the Pr\"{u}fer angles $\theta _{k}(\varphi )$
and in this case it is thus sufficient to consider 
\begin{equation*}
\varrho =\dfrac{1}{\theta _{m+1}^{\prime }(\varphi )}\equiv s\circ \theta
_{m+1}^{\prime }(\varphi )~.
\end{equation*}%
For $m<j^{\ast }$, we have%
\begin{equation}
th_{m+1}^{\prime }(\varphi )=-2t\sin \varphi +n\theta _{m+1}^{\prime
}=-2t\sin \varphi \left( 1+O(1)\right)  \label{th}
\end{equation}%
and the derivatives of higher order%
\begin{equation*}
th_{m+1}^{(k)}(\varphi )-n\theta _{m+1}^{(k)}=\left\{ 
\begin{array}{lll}
2t\left( -1\right) ^{(k+1)/2}\sin \varphi & \text{if} & k\text{ is odd} \\ 
2t\left( -1\right) ^{k/2}\cos \varphi & \text{if} & k\text{ is even}%
\end{array}%
\right.
\end{equation*}%
satisfies, in view of (\ref{ineq}),%
\begin{equation}
th_{m+1}^{(k)}(\varphi )\leq 2\left\vert t\right\vert +nC_{k}\eta ^{k}\beta
_{m}^{k}k!~.  \label{hk}
\end{equation}%
It is also sufficient to consider in both cases 
\begin{eqnarray*}
f_{0} &=&\dfrac{1}{R_{m}^{2}(\varphi )} \\
&=&\frac{1}{R_{0}^{2}}\prod_{j=1}^{m}\frac{p^{2}}{a+b\cos 2\theta _{j}+c\sin
2\theta _{j}} \\
&\equiv &\frac{1}{R_{0}^{2}}\prod_{j=1}^{m}F\circ \theta _{j}(\varphi )
\end{eqnarray*}%
where $F(\theta )$ satisfies, by direct computation,%
\begin{equation}
F^{[k]}\leq \frac{1-\delta /\zeta }{\delta }\left( \frac{\zeta }{\delta }%
\right) ^{k}  \label{Fk}
\end{equation}%
for $k\geq 0$ and some positive number $\zeta $. We also need%
\begin{equation}
s^{[k]}(x)=\frac{(-1)^{k}}{x^{k+1}}~,\ k=1,\ldots \ .  \label{rk}
\end{equation}

\begin{proposition}
\label{estimates}Let $\theta _{k}=\theta _{k}(\varphi )$, $k=0,1,\ldots $,
be the sequence of Pr\"{u}fer angles and let $\eta $, $\delta $ and $\left\{
C_{n}\right\} $ be the constants that appear in Proposition \ref{scott}.
Then, there exist positive numbers $d$ and $\hat{\eta}$, which can be
expressed in terms of the previous constants, such that (with $\varrho
=\varrho ^{\lbrack 0]}\leq d/\beta _{m+1}$) 
\begin{equation}
\varrho ^{\lbrack n]}=\left( s\circ \theta _{m+1}^{\prime }\right)
^{[n]}\leq \frac{d}{\beta _{m+1}}C_{n}\hat{\eta}^{n}\beta _{m}^{n}
\label{rhon}
\end{equation}%
as well as ($f_{0}=f_{0}^{[0]}\leq R_{0}^{-2}\delta ^{-m}$)%
\begin{equation}
f_{0}^{[n]}=\frac{1}{R_{0}^{2}}\left( \prod_{j=1}^{m}F\circ \theta
_{j}\right) ^{[n]}\leq \frac{1}{R_{0}^{2}}C_{n}(\zeta \eta )^{n}\frac{1}{%
\delta ^{m}}\beta _{m}^{n}  \label{f0n}
\end{equation}%
hold for every non--negative integer $n$, with $\zeta $ as in (\ref{Fk}).
\end{proposition}

\noindent \textit{Proof.} These inequalities are established as in
Proposition \ref{scott}, by using the Scott's formula. We begin with (\ref%
{rhon}). If $\tilde{\eta}$ is the smallest constant such that $(i-1)^{2}\eta
^{i}/i\leq \tilde{\eta}^{i}$ holds for every $i\geq 1$, by (\ref{ineq}), (%
\ref{delta-eta}) and (\ref{sf}) with $g$ and $f$ replaced by $s$ and $\theta
_{m}^{\prime }$, we have 
\begin{eqnarray*}
\varrho ^{\lbrack n]} &=&\sum_{k=1}^{n}s^{[k]}\circ \theta _{m+1}^{\prime
}\sum_{\substack{ i_{1},\ldots ,i_{k}\geq 1  \\ i_{1}+\cdots +i_{k}=n}}%
(i_{1}+1)\theta _{m+1}^{[i_{1}+1]}\cdots (i_{k}+1)\theta _{m+1}^{[i_{k}+1]}
\\
&\leq &\tilde{\eta}^{n}\beta _{m}^{n}\sum_{k=1}^{n}\frac{1}{\left( \theta
_{m+1}^{\prime }\right) ^{k+1}}\tilde{\eta}^{k}\beta _{m}^{k}\sum_{\substack{
i_{1},\ldots ,i_{k}\geq 1  \\ i_{1}+\cdots +i_{k}=n}}C_{i_{1}}\cdots
C_{i_{k}} \\
&\leq &dC_{n}\hat{\eta}^{n}\frac{\beta _{m}^{n}}{\beta _{m+1}}
\end{eqnarray*}%
where $\hat{\eta}=\delta \tilde{\eta}^{2}/(1-\Delta )$ and $d=\delta \tilde{%
\eta}/\left( \delta \tilde{\eta}+\Delta -1\right) $. In the third inequality
we have used the lower bound (\ref{beta-beta}) for $\theta _{m}^{\prime }$
and (\ref{C}). Note $\delta \tilde{\eta}>\delta \eta >1$, by definition of $%
\eta $ in Proposition \ref{scott}.

For (\ref{f0n}), we start with the Scott's formula (\ref{sf}) with $g$ and $%
f $ replaced by $F$ and $\theta _{j}$ which, together with (\ref{ineq}), (%
\ref{delta-eta}) and (\ref{Fk}), gives 
\begin{eqnarray*}
\left( F\circ \theta _{j}\right) ^{[n]} &=&\sum_{k=1}^{n}F^{[k]}\circ \theta
_{j}\sum_{\substack{ i_{1},\ldots ,i_{k}\geq 1  \\ i_{1}+\cdots +i_{k}=n}}%
\theta _{j}^{[i_{1}]}\cdots \theta _{j}^{[i_{k}]} \\
&\leq &C_{n}\left( \delta \eta \right) ^{n}\beta _{j}^{n}\frac{1-\delta
/\zeta }{\delta }\sum_{k=1}^{n}\left( \frac{\zeta }{\delta }\right) ^{k} \\
&\leq &\frac{1}{\delta }C_{n}(\zeta \eta )^{n}\beta _{j}^{n}~.
\end{eqnarray*}%
Now we take the $n$--th derivative of the product. For this, we use the
variation of Lemma \ref{A} mentioned in Remark \ref{variation}: 
\begin{eqnarray}
\left( \prod_{j=1}^{m}F\circ \theta _{j}(\varphi )\right) ^{[n]} &=&\sum 
_{\substack{ n_{1},\ldots ,n_{m}\geq 0  \\ n_{1}+\cdots +n_{m}=n}}\left(
F\circ \theta _{1}\right) ^{[n_{1}]}\cdots \left( F\circ \theta _{m}\right)
^{[n_{m}]}  \notag \\
&\leq &C_{n}(\zeta \eta )^{n}\frac{1}{\delta ^{m}}\beta _{m}^{n}  \label{Cn}
\end{eqnarray}%
concluding the proof of this proposition.

\hfill $\Box $

\begin{remark}
\label{est}An estimate of (\ref{rhon}) with $\varrho =1/th_{m+1}^{\prime
}(\varphi )$ for $m<j^{\ast }$ is analogously given by%
\begin{equation}
\varrho ^{\lbrack k]}=\left( s\circ th_{m+1}^{\prime }\right) ^{[k]}\leq 
\frac{d_{n}}{\left\vert t\right\vert }C_{k}\hat{\eta}^{k}\beta _{m}^{k}
\label{rhok}
\end{equation}%
with $d_{n}=2dn/c$ where $c=\min_{\varphi \in \text{supp}f\circ \lambda
}2\sin \varphi $. Note that, for any fixed $m$ and $t$ satisfying (\ref{t}),
the second term of the l.h.s. of (\ref{hk}) rapidly overcomes $t$. On the
other hand, an estimate for $k$ in which $t$ still dominates (\ref{hk}) is,
by the Scott's formula (\ref{sf}), much better than (\ref{rhok}):%
\begin{eqnarray}
\varrho ^{\lbrack k]} &=&\frac{1}{t}\sum_{l=1}^{k}s^{[l]}\circ
h_{m+1}^{\prime }\sum_{\substack{ i_{1},\ldots ,i_{l}\geq 1  \\ i_{1}+\cdots
+i_{l}=k}}\frac{1}{i_{1}!}h_{m+1}^{(i_{1}+1)}\cdots \frac{1}{i_{l}!}%
h_{m+1}^{(i_{l}+1)}  \notag \\
&\leq &\frac{2^{k}}{\left\vert t\right\vert }\sum_{l=1}^{k}\frac{2^{l}}{%
c^{l+1}}\sum_{\substack{ i_{1},\ldots ,i_{l}\geq 1  \\ i_{1}+\cdots +i_{l}=k 
}}\frac{1}{i_{1}!}\cdots \frac{1}{i_{l}!}  \notag \\
&\leq &\frac{2^{k}}{\left\vert t\right\vert }\frac{1}{k!}\sum_{l=1}^{k}\frac{%
2^{l}l^{k}}{c^{l+1}}~.~  \label{rhokk}
\end{eqnarray}
\end{remark}

Let us put all together. Plugging (\ref{rhon}) and (\ref{f0n}) into (\ref{f0}%
), deduced in Proposition \ref{Ln} by applying $n$ times integration by
parts $n!f_{n}=L^{n}f_{0}$ to the integrand $f_{0}$ of (\ref{I}), we arrive
at the following estimate: 
\begin{eqnarray}
n!\left\vert f_{n}\right\vert &\leq &n!\sum_{_{\substack{ k_{1},\ldots
,k_{n},p_{n}\geq 0  \\ k_{1}+\cdots +k_{n}+p_{n}=n}}}\left\vert \varrho
^{\lbrack k_{1}]}\cdots \varrho ^{\lbrack k_{n}]}f_{0}^{[p_{n}]}\right\vert 
\notag \\
&\leq &\frac{1}{R_{0}^{2}}D^{n}\frac{1}{\delta ^{m}}\left( \frac{\beta _{m}}{%
\beta _{m+1}}\right) ^{n}n!\sum_{_{\substack{ k_{1},\ldots ,k_{n},p_{n}\geq
0  \\ k_{1}+\cdots +k_{n}+p_{n}=n}}}C_{k_{1}}\cdots C_{k_{n}}C_{p_{n}} 
\notag \\
&\leq &\frac{1}{R_{0}^{2}}C_{n}D^{n}\frac{1}{\delta ^{m}}\left( \frac{\beta
_{m}}{\beta _{m+1}}\right) ^{n}n!  \label{fn}
\end{eqnarray}%
where $D=d\cdot \max \left( \hat{\eta},\zeta \eta \right) $. Estimate (\ref%
{fn}) will be used to get an upper bound for all non--resonant integrals of (%
\ref{f2drho}).

\begin{remark}
As $C_{n}=K/n^{2}$, $n\geq 1$, ($C_{0}=K$) with $K\leq 1/\left( 2+2\pi
^{2}/3\right) $ are bounded constants, (\ref{fn}) makes explicit the
dependence on the number $n$ of times that integration by parts is applied
to integral of type (\ref{I}). Explicit dependence of $n$ was not necessary
in reference \cite{KR}, since $n$ is an arbitrarily large but fixed number.
Apart this, (\ref{fn}) agrees with the estimate used on p. 522 of \cite{KR}.
\end{remark}

\section{Proof of Theorem \protect\ref{onehalf}\label{PT}}

\setcounter{equation}{0} \setcounter{theorem}{0}

Let $t$ be a fixed number. We assume $t$ positive but the negative value can
be dealt similarly. Let $f$ be a $\mathcal{C}^{\infty }$ function with
compact support in $\left( 0,2\right) $ and let $I_{f}=[\varphi _{-},\varphi
_{+}]$ be smallest closed interval that contains $\text{supp}f\circ \lambda $%
, $\lambda (\varphi )=2\cos \varphi $. Since the spectral measure is
symmetric, $d\rho (-\lambda )=d\rho (\lambda )$, we need only to consider $f$
supported in one--half of the essential spectrum. We have excluded the
origin to avoid that the curvature of $\cos ^{-1}\lambda /2$ in (\ref{tx})
vanishes (see observation right before (\ref{K1})).

For $j^{\ast }=j^{\ast }(t)$ defined by equation (\ref{t}), let $n^{\ast
}=n^{\ast }(t)$ be given by%
\begin{equation}
(n^{\ast }-1)\beta _{j^{\ast }}\leq t<n^{\ast }\beta _{j^{\ast }}~.
\label{n}
\end{equation}%
Since $\beta _{j^{\ast }}/t>1/n^{\ast }$ there are at most $n^{\ast }$
points $\varphi _{1},\ldots ,\varphi _{n^{\ast }}$ in the support of $f\circ
\lambda $ satisfying%
\begin{equation}
-\sin \varphi _{l}+\frac{l\theta _{j^{\ast }}^{\prime }}{t}=0~.  \label{sin0}
\end{equation}%
Observe that $1\leq n^{\ast }(t)\leq \beta _{j^{\ast }+1}/\beta _{j^{\ast
}}+1$, by (\ref{t}) and (\ref{n}). For the sparseness increment $\beta _{j}$
in (\ref{betaj}), we have%
\begin{equation*}
j^{\ast }(t)=\frac{\ln t}{c\ln ^{2}\ln t}\left( 1+O\left( \frac{\ln \ln \ln t%
}{\ln \ln t}\right) ~\right)
\end{equation*}%
and, consequently, the number of resonant values $n^{\ast }$ is a monotone
nondecreasing function of $t$ in each interval $(\beta _{j^{\ast }},\beta
_{j^{\ast }+1}]$. Let $L(t)$ denote the continuous interpolation of $n^{\ast
}(t)$. It follows from these observations that $L$ is a piecewise linear
function with inclination $1/\beta _{j}$ satisfying 
\begin{equation}
1<L(t)<Ee^{c\ln ^{2}\ln t}\equiv \Omega ^{2}(t)~~,\qquad \beta _{j}<t\leq
\beta _{j+1}  \label{j-nstar}
\end{equation}%
for some constant $E$, independent of $t$. The inequality (\ref{j-nstar})
will be used at the end of this section.

By (\ref{f2drho}), the Fourier--Stieltjes transform of $\rho $ can be
written as%
\begin{equation*}
\widehat{\left\vert f\right\vert ^{2}d\rho }(t)=I_{j^{\ast
}-1,0}(t)+\sum_{j=j^{\ast }-1}^{\infty }\sum_{n=1}^{\infty }\left(
~I_{j,n}(t)+~\bar{I}_{j,n}(-t)\right)
\end{equation*}%
where%
\begin{equation}
I_{j,n}(t)=\frac{1}{\pi }\int_{0}^{\pi }\left\vert f(2\cos \varphi
)\right\vert ^{2}\frac{\sin ^{2}\varphi }{R_{j}^{2}}A^{n}(\varphi
)e^{2it(\cos \varphi +n\theta _{j+1}/t)}d\varphi ~.  \label{Ikn0}
\end{equation}

We apply integration by parts to all terms of this sum not satisfying the
resonant condition (\ref{sin0}). Since the support of $f$ is compact, every
boundary term vanishes. Integration by parts may be repeated $N_{j}$ times
depending on the index $j$ of the sum. Propositions \ref{Ln} and \ref%
{estimates}, together with its combined estimate (\ref{fn}), can be used to
get an upper bound for each integral (\ref{Ikn0}) with $j\geq j^{\ast }$.
This yields 
\begin{equation}
\sum_{j=j^{\ast }}^{\infty }\sum_{n=1}^{\infty }\left( \left\vert
I_{j,n}\right\vert +\left\vert \bar{I}_{j,n}(-t)\right\vert \right) \leq
2\sum_{j=j^{\ast }}^{\infty }\left( \sum_{n=1}^{\infty }\frac{1}{n}%
a^{n}\right) \frac{1}{R_{0}^{2}}C_{N_{j}}D^{N_{j}}\frac{1}{\delta ^{j}}%
\left( \frac{\beta _{j}}{\beta _{j+1}}\right) ^{N_{j}}N_{j}!  \label{Ikn}
\end{equation}%
where $a=\sup_{\varphi \in I_{f}}\left\vert A(\varphi )\right\vert <1$. If
the sequences $\left( \beta _{j}\right) _{j}$ and $\left( N_{j}\right) _{j}$
are chosen so that%
\begin{equation}
D^{N_{j}}\frac{1}{\delta ^{j}}\left( \frac{\beta _{j}}{\beta _{j+1}}\right)
^{N_{j}}N_{j}!\leq \frac{1}{\beta _{j+1}}  \label{if}
\end{equation}%
then the series in (\ref{Ikn}) converges uniformly in $t$. By (\ref{delta})
and (\ref{t}), we have 
\begin{equation*}
\frac{1}{\beta _{j+1}}=\frac{1}{\beta _{j^{\ast }+1}}\frac{\beta _{j^{\ast
}+1}}{\beta _{j^{\ast }+2}}\cdots \frac{\beta _{j}}{\beta _{j+1}}\leq \delta
^{j-j^{\ast }}\frac{1}{\beta _{j^{\ast }+1}}\ .
\end{equation*}%
and $t/\beta _{j^{\ast }+1}<1$. Consequently, 
\begin{eqnarray}
\sum_{j=j^{\ast }}^{\infty }\sum_{n=1}^{\infty }\left( \left\vert
I_{j,n}\right\vert +\left\vert \bar{I}_{j,n}(-t)\right\vert \right)  &\leq &%
\frac{2a}{1-a}\frac{1}{R_{0}^{2}}\sum_{j=j^{\ast }}^{\infty }C_{N_{j}}\frac{1%
}{\beta _{j+1}}  \notag \\
&\leq &\frac{2a}{1-a}\frac{1}{R_{0}^{2}}\frac{t}{\beta _{j^{\ast }+1}}\frac{K%
}{t}\sum_{l=0}^{\infty }\delta ^{l}\leq \frac{C}{t}  \label{C0}
\end{eqnarray}%
holds with $C<\infty $ independent of $t$.

Let us now verify that the sparseness condition stated in Theorem \ref%
{onehalf} satisfies (\ref{if}). Choosing $\beta _{j}$ as given by (\ref%
{betaj}) and $N_{j}=j+1$, by the Stirling formula,%
\begin{equation*}
\beta _{j+1}\left( \frac{\beta _{j}}{\beta _{j+1}}\right) ^{N_{j}}N_{j}!=%
\frac{\delta ^{j}}{\sqrt{2\pi (j+1)}}\left( \frac{j+1}{e}\right) ^{j+1}\exp
\left( -2c(j+1)\ln \left( j+1\right) \right) \left( 1+O\left( \frac{\ln j}{j}%
\right) \right)
\end{equation*}%
and (\ref{if}) holds for any $c>1/2$ provided $j$ is large enough.

Note that, for the sparseness increment $\left( \beta _{j}\right) _{j\geq 1}$
given by (\ref{factorial}),%
\begin{equation*}
\beta _{j+1}^{1-\varepsilon }\left( \frac{\beta _{j}}{\beta _{j+1}}\right)
^{N_{j}}=\delta ^{\varepsilon (j+1)}\exp \left( -(j+1)\ln \left( j+1\right)
-\varepsilon ^{-1}(j+1)\right) <\frac{1}{N_{j}!}\delta ^{j}D^{-N_{j}}
\end{equation*}%
and (\ref{if}) holds with $1/\beta _{j+1}$ replaced by $1/\beta
_{j+1}^{1-\varepsilon }$, provided $1+1/\varepsilon >\ln D-(1-\varepsilon
)\ln \delta $ and $j$ is large enough. Together with Remark \ref{tepsilon},
the proof may continued exactly as for decaying $t^{-1}$. We shall consider
only the latter case.

For $j<j^{\ast }$, (\ref{th}) holds and we need replace the estimate (\ref%
{rhon}) by (\ref{rhok}) and $t$ occupies now the place of $\beta _{m+1}$ in (%
\ref{fn}). Applying successive integration by parts to $I_{j^{\ast }-1,0}$
gives, analogously%
\begin{equation}
\left\vert I_{j^{\ast }-1,0}\right\vert \leq \frac{1}{R_{0}^{2}}%
C_{N_{j^{\ast }-1}}D^{N_{j^{\ast }-1}}\frac{1}{\delta ^{j^{\ast }-1}}\left( 
\frac{\beta _{j^{\ast }-1}}{t}\right) ^{N_{j^{\ast }-1}}N_{j^{\ast }-1}!\leq 
\frac{C^{\prime }}{t}  \label{I0}
\end{equation}%
for some constant $C^{\prime }$. Note that, by (\ref{t}),%
\begin{equation*}
\left( \frac{\beta _{j^{\ast }-1}}{t}\right) ^{N_{j^{\ast }-1}}\leq \frac{1}{%
t}\frac{\beta _{j^{\ast }-1}^{N_{j^{\ast }-1}}}{\beta _{j^{\ast
}}^{N_{j^{\ast }-1}-1}}
\end{equation*}%
and by (\ref{if})%
\begin{equation}
D^{N_{k}}\frac{1}{\delta ^{k}}\frac{\beta _{k}^{N_{k}}}{\beta
_{k+1}^{N_{k}-1}}N_{k}!\leq 1~  \label{if2}
\end{equation}%
for $k$ large enough.

It remains to estimate the sum $S_{j^{\ast }}(t)=\displaystyle%
\sum_{n=1}^{\infty }~\left( I_{j^{\ast }-1,n}(t)+\bar{I}_{j^{\ast
}-1,n}(-t)\right) $ which contains the most significant terms responsible
for $t^{-1/2}$ decaying behavior. To extract this decay we write%
\begin{equation}
I_{j^{\ast }-1,n}(t)=\frac{1}{\pi }\int_{0}^{\pi }\left\vert f(2\cos \varphi
)\right\vert ^{2}\frac{\sin ^{2}\varphi }{R_{j^{\ast }-1}^{2}}B_{j^{\ast
}-1}^{n}(\varphi )e^{2it(\cos \varphi +n\beta _{j^{\ast }}\varphi
/t)}d\varphi  \label{Ijstar}
\end{equation}%
(analogously for $\bar{I}_{j^{\ast }-1,n}(t)$) where, by (\ref{thetatheta}), 
$B_{k}$ is a function of the Pr\"{u}fer angles $\theta _{k}(\varphi )$ such
that $\left\vert B_{k}\right\vert =\left\vert A\right\vert $ and%
\begin{eqnarray*}
\arg B_{k} &=&\arg A+\theta _{k+1}(\varphi )-\beta _{k+1}\varphi \\
&=&\arg A+g\circ \theta _{k}(\varphi )
\end{eqnarray*}%
with $g$ given by (\ref{gg}). We then apply Lemma \ref{main} to (\ref{Ijstar}%
) with 
\begin{equation*}
d\left( G_{n}\circ \lambda \right) (\varphi )=\frac{1}{\pi }\left\vert
f(2\cos \varphi )\right\vert ^{2}\frac{\sin ^{2}\varphi }{R_{j^{\ast }-1}^{2}%
}B_{j^{\ast }-1}^{n}(\varphi )d\varphi
\end{equation*}%
and%
\begin{equation*}
tx(t,\lambda )=t\lambda +2n\beta _{j^{\ast }}\cos ^{-1}\frac{\lambda }{2}~.
\end{equation*}%
Note that, by (\ref{t}), $\kappa =2\pi n\beta _{j^{\ast }}$ and $n^{\ast
}\beta _{j^{\ast }}=O(t)$. In order to fulfill all assumptions of Lemma \ref%
{main} it remains to show that $\widehat{dG_{n}}(t)$ decays as $\left\vert
t\right\vert ^{-1}$ (see equation (\ref{Ghat})).

We estimate the Fourier--Stieltjes transform 
\begin{equation}
\widehat{dG_{n}}(t)=\frac{1}{\pi }\int_{0}^{\pi }\left\vert f(2\cos \varphi
)\right\vert ^{2}\frac{\sin ^{2}\varphi }{R_{j^{\ast }-1}^{2}}B_{j^{\ast
}-1}^{n}(\varphi )e^{2it\cos \varphi }d\varphi  \label{dG}
\end{equation}%
as the non--resonant integrals (\ref{Ikn0}) with $j<j^{\ast }$ (see Remark %
\ref{est}). The estimate (\ref{rhon}) is replaced by (\ref{rhokk}) and $%
f_{0}=\left\vert f(2\cos \varphi )\right\vert ^{2}\dfrac{\sin ^{2}\varphi }{%
R_{j^{\ast }-1}^{2}(\varphi )}A^{n}(\varphi )\exp \left( ing\circ \theta
_{j^{\ast }-1}(\varphi )\right) $ includes now an extra exponential term
depending on $\theta _{j^{\ast }-1}(\varphi )$.

To deal with this new term we need some more estimates. By the Scott's
formula (\ref{sf}) together with (\ref{g}), we have%
\begin{eqnarray}
\left\vert \left( e^{ing}\right) ^{[N]}\right\vert &\leq &\sum_{k=1}^{N}%
\frac{n^{k}}{k!}\sum_{\substack{ i_{1},\ldots ,i_{k}\geq 1  \\ i_{1}+\cdots
+i_{k}=N}}\left\vert g^{[i_{1}]}\cdots g^{[i_{k}]}\right\vert  \notag \\
&\leq &c_{2}\xi ^{N}  \label{c2}
\end{eqnarray}%
with $c_{2}=e^{nc_{1}}-1$. The Scott's formula (\ref{sf}) applied once again
together with Proposition \ref{scott}, (\ref{delta-eta}) and (\ref{c2}),
yield%
\begin{eqnarray*}
\left\vert \exp \left( ing\circ \theta _{j^{\ast }-1}\right)
^{[N]}\right\vert &\leq &\sum_{k=1}^{N}\left\vert \left( e^{ing}\right)
^{[k]}\circ \theta _{j^{\ast }-1}\right\vert \sum_{\substack{ i_{1},\ldots
,i_{k}\geq 1  \\ i_{1}+\cdots +i_{k}=N}}\left\vert \theta
_{j}^{[i_{1}]}\right\vert \cdots \left\vert \theta _{j}^{[i_{k}]}\right\vert
\\
&\leq &c_{3}C_{N}\left( \delta \xi \eta \right) ^{N}\beta _{j^{\ast
}-1}^{N}~,
\end{eqnarray*}%
with $c_{3}=c_{2}\xi /(\xi -1)$. Finally, we shall replace (\ref{Cn}) by 
\begin{eqnarray*}
\left\vert \left( \prod_{j=1}^{j^{\ast }-1}F\circ \theta _{j}\cdot
e^{ing\circ \theta _{j^{\ast }-1}}\right) ^{[N]}\right\vert &\leq &\sum 
_{\substack{ n_{1},\ldots ,n_{j^{\ast }}\geq 0  \\ n_{1}+\cdots +n_{j^{\ast
}}=N}}\left\vert \left( F\circ \theta _{1}\right) ^{[n_{1}]}\cdots \left(
F\circ \theta _{j^{\ast }-1}\right) ^{[n_{j^{\ast }-1}]}\left( e^{ing\circ
\theta _{j^{\ast }-1}}\right) ^{[n_{j^{\ast }}]}\right\vert \\
&\leq &c_{3}C_{N}(\bar{\zeta}\eta )^{N}\frac{1}{\delta ^{j^{\ast }-1}}\beta
_{j^{\ast }-1}^{N}
\end{eqnarray*}%
with $\bar{\zeta}=\max \left( \delta \xi ,\zeta \right) $.

Since the require estimates didn't change significantly, we integrate (\ref%
{dG}) by parts $N_{j^{\ast }-1}$ times and use (\ref{fn}) with (\ref{rhokk})
in the place of (\ref{rhon}) to get, exactly as for (\ref{I0}),%
\begin{equation*}
\left\vert \widehat{dG_{n}}(t)\right\vert \leq \frac{C^{\prime \prime }}{t}%
d^{k}\tilde{a}^{n}\ ,\qquad \beta _{k}\leq t<\beta _{k+1}
\end{equation*}%
for some constants $C^{\prime \prime }<\infty $, $d<1$ and $\tilde{a}%
=e^{c_{1}}\sup_{\varphi \in I_{f}}\left\vert A(\varphi )\right\vert <1$, as $%
c_{1}$ is arbitrarily small by the observation after (\ref{g}). Here, we
have used the fact that (\ref{if2}) holds with $D$ replaced by $D/d$, for
any $d>1$, provided $k$ is large enough. This immediately imply, by a slight
modification of Lemma \ref{main} (see Remark \ref{tepsilon}),%
\begin{eqnarray*}
\left\vert I_{j^{\ast }-1,n}(t)\right\vert  &\leq &\int_{\left\vert \tau
\right\vert \leq \Delta \left\vert t\right\vert }\left\vert \Lambda (t,\tau
)\right\vert \left\vert \widehat{dG}(\tau )\right\vert d\tau +O\left(
1/t\right)  \\
&\leq &\frac{2K}{\sqrt{\kappa }}\tilde{a}^{n}\sum_{k=0}^{j^{\ast
}}d^{k}\int_{\beta _{k}}^{\beta _{k+1}}\frac{C^{\prime \prime }}{\tau }d\tau
+O\left( 1/t\right)  \\
&\leq &\frac{2KC^{\prime \prime }}{\sqrt{n^{\ast }\kappa }}n^{\ast }\tilde{a}%
^{n}\sum_{k=0}^{\infty }d^{k}\ln \frac{\beta _{k+1}}{\beta _{k}}+O\left(
1/t\right) 
\end{eqnarray*}%
and by (\ref{n}) and the fact that $\ln \beta _{k+1}/\beta _{k}=O(\ln ^{2}k)$%
, we have 
\begin{eqnarray*}
\left\vert S_{j^{\ast }}(t)\right\vert  &\leq &\sum_{n=1}^{\infty }~\left(
\left\vert I_{j^{\ast }-1,n}(t)\right\vert +\left\vert \bar{I}_{j^{\ast
}-1,n}(-t)\right\vert \right)  \\
&\leq &\frac{C}{\sqrt{\left\vert t\right\vert }}\Omega (\left\vert
t\right\vert )~,
\end{eqnarray*}%
where $\Omega (t)$ is defined in (\ref{j-nstar}). 

Now, we show that $\Omega (t)$ increases slower than $t^{\varepsilon }$, for
any $\varepsilon >0$. Suppose, by contradiction, that $\lim_{t\rightarrow
\infty }\Omega (t)/t^{\varepsilon }=k>0$ holds for some $\varepsilon >0$.
Then, by L'Hospital,%
\begin{equation*}
\lim_{t\rightarrow \infty }\frac{\Omega (t)}{t^{\varepsilon }}%
=\lim_{t\rightarrow \infty }\frac{\Omega ^{\prime }(t)}{\varepsilon
t^{\varepsilon -1}}=\frac{c}{\varepsilon }\lim_{t\rightarrow \infty }\frac{%
\Omega (t)}{t^{\varepsilon }}\cdot \lim_{t\rightarrow \infty }\frac{\ln \ln t%
}{\ln t}=0~,
\end{equation*}%
concluding the proof of Theorem \ref{onehalf}.

\hfill $\Box $

\end{document}